\documentclass[12pt,a4paper]{article}
\usepackage{amsmath,amssymb,amsthm}
\usepackage{fullpage}

\newtheorem{thm}{Theorem}

\newcommand{\dis}{\displaystyle}
\newtheorem{example}{Example}
\newtheorem{lemma}{Grossehilfsatz}
\numberwithin{equation}{subsection}
\begin{document}

\title{Mertens' Proof of Mertens' Theorem}
\author{Mark B. Villarino\\
Depto.\ de Matem\'atica, Universidad de Costa Rica,\\
2060 San Jos\'e, Costa Rica}
\date{April 28, 2005}

\maketitle

 \begin{abstract}
 We study\textsc{ Mertens}' own proof  (1874) of his theorem on the sum of the reciprocals of the primes and compare it with the modern treatments.
 \end{abstract}
\tableofcontents
\section{Historical Introduction}
\subsection{Euler}
In 1737, \textsc{Leonhard Euler} created analytic (prime) number theory with the publication of
his memoir ``Variae observationes circa series infinitas" in \textit{Commentarii academiae scientiarum Petropolitanae} 9 (1737), 160-188; \textit{Opera omnia} (1) XIV, 216-244.  \textbf{Theorema 7} states:
\begin{quote}
\emph{``If we take to infinity the continuation of these fractions
$$\frac{2\cdot3\cdot5\cdot7\cdot11\cdot13\cdot17\cdot19\cdots}{1\cdot2\cdot4\cdot6\cdot10\cdot12\cdot16\cdot18\cdots}
$$
where the numerators are all the prime numbers and the denominators are the numerators less one unit, the result is the same as the sum of the series
$$1+\frac{1}{2}+\frac{1}{3}+\frac{1}{4}+\frac{1}{5}+\frac{1}{6}+\cdots."
$$}
\end{quote}

This is the wonderful identity which, today, we write \cite{Hardy}, \cite{Ingham}, \cite{Jameson}:
\begin{align}
\fbox{$\dis \prod_{2}^{\infty}\frac{1}{1-\dfrac{1}{p^{1+\rho}}}=1+\frac{1}{2^{1+\rho}}+\frac{1}{3^{1+\rho}}+\frac{1}{4^{1+\rho}}+\cdots,$}
\end{align}Here $\rho>0$ and the product on the left is taken over all primes $p\geqslant 2$, while the right hand side is the famous \textsc{Riemann} zeta function, $\zeta(1+\rho).$  The modern statement is nice, but does not have the sense of wonder that \textsc{Euler}'s statement carries.  Yes, it is not rigorous, but it \emph{is} beautiful.

\textsc{Euler}'s memoir is replete with extraordinary identities relating infinite products and series of primes, but our interest is in his \textbf{Theorema 19}: 
\begin{quote}
\emph{``Summa seriei reciprocae numerorum primorum
$$\frac{1}{2}+\frac{1}{3}+\frac{1}{5}+\frac{1}{7}+\frac{1}{11}+\frac{1}{13}+\rm etc.
$$
est infinite magna, infinities tamen minor quam summa seriei harmonicae
$$1+\frac{1}{2}+\frac{1}{3}+\frac{1}{4}+\frac{1}{5}+\rm etc.
$$
Atque illius summa est huius summae quasi logarithmus."}
\end{quote}

We translate this as our first formal theorem.
\begin{thm}
The sum of the reciprocals of the prime numbers
$$\frac{1}{2}+\frac{1}{3}+\frac{1}{5}+\frac{1}{7}+\frac{1}{11}+\frac{1}{13}+\rm etc.
$$
is infinitely great but is infinitely times less than the sum of the harmonic series
$$1+\frac{1}{2}+\frac{1}{3}+\frac{1}{4}+\frac{1}{5}+\rm etc.
$$
And the sum of the former is as the logarithm of the sum of the latter.
\end{thm}\hfill$\Box$

The last line of \textsc{Euler}'s attempted proof is:
\begin{quote}
\emph{``. . . and finally,
$$\frac{1}{2}+\frac{1}{3}+\frac{1}{5}+\frac{1}{7}+\frac{1}{11}+\cdots=\ln\ln \infty".
$$}
\end{quote}
(We have written ``$\ln\ln \infty$" instead of \textsc{Euler}'s ``\emph{ll}$\infty$.")

It is evident that \textsc{Euler} says that the series of prime reciprocals \emph{diverges} and that the partial sums grow like the \emph{logarithm} of the partial sums of the harmonic series, that is $\sum_{p\leqslant x}\frac{1}{p}$ grows like $\ln\ln x.$  Of course, this implies (trivially) that \emph{there are infinitely many primes}, since the series of reciprocal primes must necessarily have infinitely many summands.  Moreover, it even indicates \emph{the velocity} of divergence and therefore \emph{the density} of the primes, a totally new idea.

This was the first application of analysis (limits and infinite series) to prove a theorem in number theory, the first \emph{new} proof of the infinity of primes in two thousand years (!), and opened an entirely new branch of mathematics, analytic number theory, which is a rich and fecund area of modern mathematics.
\subsection{Legendre and Chebyshev}
The first quantitative statement of \textsc{Euler}'s theorem on the sum of the reciprocal primes appeared in \textsc{Legendre}'s \emph{Th\'eorie des nombres} (troisi\`eme \'edition, quatri\`eme partie, VIII, (1808)), namely:
$$\sum_{p\leqslant G}\frac{1}{p}= \ln(\ln G-0.08366)+C,
$$
where $G$ is a given real number and $C$ is an unknown numerical constant.  \textsc{Legendre} gave no hint of a proof nor of the origin of the mysterious constant ``0.08366."

In 1852, no less a mathematician than the great russian analyst \textsc{Chebyshev} \cite{Chebyshev} attempted a proof of \textsc{Legendre}'s theorem, but failed.  The problem of finding such a proof became celebrated, and the stage was set for its solution.
\subsection{Mertens}
In 1874 (see \cite{Mertens}) the brilliant young Polish-Austrian mathematician \footnote{He was a professor of mathematics for over 20 years (1865-1884)  at the Jagiellonian university in Cracow.  At that time, Poland was partitioned among Prussia, Russia and Austria, and Cracow was in the austrian zone -- there was not an independent polish state then.  Mertens' wife was polish and he spoke polish as well as german. Then he went to Graz to become rector of the politechnique there. \cite{Wolf}}, \textsc{Franciszek Mertens}, published a proof of his now famous theorem on the sum of the prime reciprocals:

\begin{thm}(Mertens (1874))
Let $x\geqslant 1$ be any real number.  Then
\begin{align}
\fbox{$\dis\sum_{p\leqslant x}\frac{1}{p}= \ln\ln[x]+\gamma+\sum_{m=2}^\infty \mu(m)\frac{\ln\{\zeta(m)\}}{m}+\delta $}
\end{align}
where $\gamma$ is \textsc{Euler}'s constant, $\mu(m)$ is the \textsc{M\"obius} function, $\zeta(m)$ is the \textsc{Riemann} zeta function, and 
\begin{equation}
\fbox{$\dis |\delta|< \frac{4}{\ln([x]+1)}+\frac{2}{[x]\ln[x]}.$}
\end{equation}
\end{thm}$\hfill\Box$

(We write $[x]:=$ the greatest integer in $x.$) We have slightly altered his notation. 

Today we write the statement of \textsc{Mertens}' theorem in the form \cite{Hardy}, \cite{Ingham}: 
\begin{thm}
$$ B:=\lim_{x\rightarrow \infty} \left(\sum_{p\leqslant x}\frac{1}{p}-\ln\ln x \right)
$$
is a well-defined constant.$\hfill\Box$
\end{thm}

An alternative more precise statement of the modern theorem is:

\begin{thm}
$$ \sum_{p\leqslant x}\frac{1}{p}=\ln\ln x+B+O\left(\frac{1}{\ln x}\right)
$$
where
$$B:=\sum_{p}\left\{\ln \left(1-\frac{1}{p}\right)+\frac{1}{p}\right\}.
$$\hfill$\Box$
\end{thm}

The modern presentations of \textsc{Mertens}' theorem, \cite{Hardy},\cite{Ingham}, 
 \cite{Jameson}
\cite{Landau},  include:\begin{enumerate}
  \item \emph{no} discussion of an explicit numerical error estimate (such as \textsc{Mertens}' $\delta$).
  \item \emph{no} computation of $B$, in particular, a proof of the wonderful formula: 
  \begin{align}
  B= \gamma+\sum_{n=2}^\infty \mu(n)\frac{\ln\{\zeta(n)\}}{n}.
  \end{align}
 \textsc{Mertens} used this formula to compute the value:
$$B\approx 0.2614972128.
$$

  \item \emph{no} hint of how \textsc{Mertens}, himself, proved his explicit theorem.
\end{enumerate}

In this paper we will present a self-contained motivated exposition of \textsc{Mertens}' original proof and compare its strategy, tactics, and details with the modern approach.  \textsc{Mertens}' proof is brilliant, insightful, and instructive.  It deserves to be better known and our paper attempts to achieve this.\footnote{\textsc{Mertens}' paper also contains a proof of his (almost) equally famous product-theorem: $$\prod_{p\leqslant G}\frac{1}{1-\frac{1}{p}}=e^{\gamma+\delta'}\cdot \ln G$$ where $|\delta'|<\frac{4}{\ln (G+1)}+\frac{2}{G\ln G}+\frac{1}{2G}$. But there is nothing new in his treatment that does not appear in the theorem we are dealing with, so we do not discuss it here.}
\section{The Modern Proof}
\subsection{Partial Summation}

Modern prime number theory, indeed number theory in general, has developed a systematic approach to the computation of finite sums of number theoretic functions by use of what is called ``Abel summation," or ``partial summation."  We follow \cite{Jameson}.

\begin{thm}(Abel Summation)
Let $y<x,$ and let $f$ be a function (with real or complex values) having a continuous derivative on $[y,x]$.  Then

\begin{align}
\sum_{y<r\leqslant x}a(r)f(r)=A(x)f(x)-A(y)f(y) -\int_{y}^x A(t)f'(t)~dt
\end{align}
where the integers $a(r)$ are given, and where
\begin{align}
A(x):=\sum_{r\leqslant x}a(r).
\end{align}

\end{thm}$\hfill\Box$

We will apply this technique to the sum 
\begin{align}
\sum_{p\leqslant x}\frac{1}{p}.
\end{align}
\subsection{The Relation with $\pi(x)$}
\begin{example} \rm Take $y:=2,$ in the Partial Summation formula, and take

\[a(r):= \left\{ \begin{array}{r@{\quad if \quad}l}
1 & r=p\\0 & r\neq p
\end{array}\right.\]
that is, $a(r)$ is the \emph{characteristic function} of the prime numbers $p$.  Moreover, take:$$f(r):=\frac{1}{r}.$$  Then we conclude that $$A(x)$$ is equal to the number of prime numbers $p\leqslant x$, i.e., the  prime counting function $\pi(x)$.  Therefore, the formula for Abel summation gives us:
\begin{align}
\sum_{p\leqslant x}\frac{1}{p}=\frac{\pi(x)}{x}+\int_{2}^{x}\frac{\pi(t)}{t^{2}}~dt,
\end{align}
a very pretty equation relating our sum to the famous function $\pi(x)$.  Unfortunately, in order to apply it we must know upper and lower bounds for $\pi(x)$, and the study of such bounds is the subject of the Prime Number Theorem, something much deeper than our topic.

\end{example}
\subsection{The First Grossehilfsatz}
\begin{example}\rm Following \textsc{Mertens} (in a slightly different context: see 3.4) we   again take $y:=2$, but this time we take  

\[a(r):= \left\{ \begin{array}{r@{\quad if \quad}l}
\dfrac{\ln p}{p} & r=p\\0 & r\neq p
\end{array}\right.\] and $$f(r):=\frac{1}{\ln r}.$$  Then
\begin{align}
A(x)=\sum_{p\leqslant x}\frac{\ln p}{p}.
\end{align}
Therefore, the formula for Abel summation gives us:
\begin{align}
\sum_{p\leqslant x}\frac{1}{p}=\frac{A(x)}{\ln x}+\int_{2}^{x}\frac{A(t)}{t(\ln t)^{2}}~dt,
\end{align} a nice formula, but with $A(x)$ the slightly more exotic function given in (2.3.1).  In his paper, \textsc{Mertens} proves two ``\textbf{\emph{Grossehilfs\"atze}}" (in \textsc{Landau}'s marvelous German phraseology: the English ``fundamental lemmas" does not carry the same force.)  The first one deals with our $A(x).$
\begin{lemma}\begin{align}
\fbox{$\dis\sum_{p\leqslant x}\frac{\ln p}{p}=\ln x +R(x), \rm where\  |R(x)|<2.$}
\end{align}$\hfill\Box$
\end{lemma}

The interest in this is the explicit numerical error estimate, $|R(x)|<2,$ which, as we will see, is quite good.

We will give \textsc{Merten}'s nice proof of this result later on (see 3.5), but for now we assume it to be true.  

Then, if we put ,$$R(t):=\sum_{p\leqslant x}\frac{\ln p}{p}-\ln t,$$ which means (by (2.3.3)) that $$|R(t)|<2,$$by (2.3.2), we conclude that 
\begin{align*}
 \sum_{p\leqslant x}\frac{1}{p}& = \frac{\ln x + R(x)}{\ln x}+\int_{2}^{x}\frac{\ln t +R(t)}{t(\ln t)^{2}}~dt  \\
 & =  1+\frac{R(x)}{\ln x}+\int_{2}^x \frac{1}{t\ln t}~dt+\int_{2}^{x}\frac{R(t)}{t(\ln t)^{2}}~dt\\
&= 1+\frac{R(x)}{\ln x}+\ln\ln x-\ln\ln 2+\int_{2}^\infty \frac{R(t)}{t(\ln t)^{2}}~dt-\int_{x}^\infty \frac{R(t)}{t(\ln t)^{2}}~dt\\
&= \ln\ln x+ \underbrace{1-\ln\ln 2+\int_{2}^\infty \frac{R(t)}{t(\ln t)^{2}}~dt}_{\textstyle{\rm a \rm \ constant  \ B}}+\underbrace{\frac{R(x)}{\ln x}-\int_{x}^\infty \frac{R(t)}{t(\ln t)^{2}}~dt}_{\textstyle{\leqslant \frac{2}{\ln x}+\frac{2}{\ln x}=\frac{4}{\ln x}}}\\
&= \ln\ln x + B +\delta,
\end{align*}where $$|\delta|<\frac{4}{\ln x}.$$We have proved:  

\begin{thm}
There exists a constant, $B$, such that for all real numbers $x\geqslant 2$, 
\begin{align}
\sum_{p\leqslant x}\frac{1}{p}=\ln\ln x + B +\delta,
\end{align}where
\begin{align}
|\delta|<\frac{4}{\ln x}.
\end{align}        \hfill$\Box$      \end{thm}
This is an explicit form of \textsc{Mertens}' theorem (our \textbf{Theorem 2}) with a somewhat \emph{better} error term than (1.2) in \textsc{Mertens}' original statement! Unfortunately, the form of the constant $$B:=1-\ln\ln 2+\int_{2}^\infty \frac{R(t)}{t(\ln t)^{2}}~dt$$ gives no clue as to how to compute it, much less that it has the form $\gamma+C$, for some constant $C$, as we saw in equation (1.3.3). This shows both the advantage, and the disadvantage of the modern approach:  it is systematic and gives a (slightly) better error term with little effort, but it gives no algorithm for the explicit computation of the constant $B.$

There are modern treatments \cite{Hardy}, \cite{Ingham}, \cite{Landau}, that show the formula $$B=\gamma+C$$ to be valid, but there  is no modern textbook treatment of the formula (1..3).  There is a beautiful recent paper \cite{L-P} on this formula and its computation which should be consulted.

\end{example}
\section{Mertens' Proof}
\subsection{A Sketch of the Proof}
\textsc{Mertens} starts with the \emph{convergent} ``prime zeta function" 
$$\sum_{p} \frac{1}{p^{1+\rho}}
$$
where $\rho>0,$ and writes its \emph{partial sum} for primes $p\leqslant x$ as:
\begin{align}
\sum_{p\leqslant x}\frac{1}{p^{1+\rho}}=\sum_{p}\frac{1}{p^{1+\rho}}-\sum_{p> x}\frac{1}{p^{1+\rho}}
\end{align}
and then studies the RHS as $\rho \rightarrow 0.$
It is fairly easy to show that 
\begin{align}
\sum_{p}\frac{1}{p^{1+\rho}}=\ln\left(\frac{1}{\rho}\right)-H+o(\rho),
\end{align}
where 
\begin{align}
H:=\sum_{n=2}^\infty \mu(n)\frac{\ln\{\zeta(n)\}}{n}
\end{align}

It \emph{takes work(!)} to show that the ``remainder,"
\begin{align}
\sum_{p> x}\frac{1}{p^{1+\rho}}= \ln\left(\frac{1}{\rho}\right)-\ln\ln x-\gamma+\delta+o(\rho).
\end{align}

Equations $(3.1.1), \ (3.1.2)$, and $(3.14)$ show 
\begin{align}
\sum_{p\leqslant x}\frac{1}{p^{1+\rho}}=\ln\ln x +\gamma-H+\delta+o(\rho),
\end{align}
and letting $\rho \rightarrow 0$ gives \textsc{Mertens}'s theorem.

The equations $(3.1.2)$ and $(3.1.4)$ show that the ``Mertens constant," $B$, is the sum of \emph{two} constants, $\gamma$ and $-H$, and each comes from a \emph{different part} of the ``prime zeta function."  It is this fact that makes \textsc{Mertens}' theorem hard to prove.

Our presentation follows \textsc{Mertens} quite closely, although we fill in several details.  His mathematics is striking and beautiful, a \emph{tour de force} of classical analysis.

\subsection{Euler-Maclurin and Stirling}
In this section we will cite the versions of the Euler-Maclaurin formula and Stirling's formula which will be used in \textsc{Mertens}'s proof. The proof of both can be found in \cite{Knopp}.
\begin{thm}(Euler-Maclaurin)
Let $f(t)$ have a continuous derivative, $f'(t)$, for $t\geqslant 1$.  Then:
\begin{align}
\sum_{n\leqslant x}f(n)=\int_{1}^{x}f(t)~dt+\int_{1}^{x}(t-[t])f'(t)~dt+f(1)-(x-[x])f(x).
\end{align}
\end{thm}\hfill$\Box$

\begin{thm}(Stirling's Formula)
The following relations are valid for all real $x\geqslant 4$ and all integers $n\geqslant 5$:
\begin{align}
\ln(1\cdot2\cdot3\cdots[x]) & <  x\ln x+\frac{1}{2}\ln x-x+\ln \sqrt{2\pi}+\frac{1}{12x} \\
2\ln\left(1\cdot2\cdot3\cdots\left[\frac{x}{2}\right]\right) & >  x\ln x-x\ln 2-\ln x-x+d\ln \sqrt{2\pi}+\ln2 -\frac{2}{x-2}\\
\ln(n!)&= n\ln n-n+\frac{1}{2}\ln n+\ln \sqrt{2\pi}+\frac{\lambda}{12n},  \ |\lambda|<1 
\end{align}
\end{thm}\hfill$\Box$


\subsection{The First Step of \textsc{Mertens'} Proof}

\textsc{Mertens} begins with \textsc{Euler}'s marvelous identity:
\begin{align}
\fbox{$\dis \prod_{2}^{\infty}\frac{1}{1-\dfrac{1}{p^{1+\rho}}}=1+\frac{1}{2^{1+\rho}}+\frac{1}{3^{1+\rho}}+\frac{1}{4^{1+\rho}}+\cdots,$}
\end{align} as indeed does most of analytic prime number theory.  Here $\rho>0$ and the product on the left is taken over all primes $p\geqslant 2.$  The right hand side is the famous \textsc{Riemann} zeta function $\zeta(1+\rho)$.

Now,
\begin{align*}
\zeta(1+\rho) & :=  \sum_{n\geqslant 1}\frac{1}{n^{1+\rho}} \\
 & \overset{(3.2.1)}{=}  \int_{1}^{\infty}\frac{1}{x^{1+\rho}}~dx+1+\theta(-1), \ \theta\in[0,1]\\
 &= -\frac{1}{\rho x^{\rho}}\bigg|_{x=1}^{x=\infty}+1-\theta\\
 &= \frac{1}{\rho}+1-\theta\\
 &= \frac{1+o(\rho)}{\rho},
\end{align*} thus
\begin{align}
\prod_{2}^{\infty}\frac{1}{1-\dfrac{1}{p^{1+\rho}}}=\frac{1+o(\rho)}{\rho}.
\end{align}

Taking logarithms of both sides we obtain:
\begin{align*}
\sum_{2}^{\infty}\ln\left(\frac{1}{1-\dfrac{1}{p^{1+\rho}}}\right) & = \sum_{2}^{\infty}\left(\frac{1}{p^{1+\rho}}+\frac{1}{2}\cdot\frac{1}{p^{2+2\rho}}+\frac{1}{3}\cdot\frac{1}{p^{3+3\rho}}+\cdots \right) \\
&=\sum_{2}^{\infty}\frac{1}{p^{1+\rho}}+\frac{1}{2}\cdot\sum_{2}^{\infty}\frac{1}{p^{2+2\rho}}+\frac{1}{3}\cdot\sum_{2}^{\infty}\frac{1}{p^{3+3\rho}}+\cdots\\
 & =  \ln \left\{\frac{1+o(\rho)}{\rho}\right\} 
\end{align*} Therefore,\begin{align}
\sum_{2}^{\infty}\frac{1}{p^{1+\rho}}=\ln \left\{\frac{1+o(\rho)}{\rho}\right\}-\frac{1}{2}\cdot\sum_{2}^{\infty}\frac{1}{p^{2+2\rho}}-\frac{1}{3}\cdot\sum_{2}^{\infty}\frac{1}{p^{3+3\rho}}-\cdots
\end{align}
\textsc{Mertens} wants to let $\rho\rightarrow 0$ on both sides of (3.3.3).  That way, formally, the left hand side becomes $$\sum_{p}\frac{1}{p},$$ the sum he wishes to study, while the right hand side becomes $$\lim_{\rho\rightarrow 0}\ln\left\{\frac{1+o(\rho)}{\rho}\right\}-\frac{1}{2}\cdot\sum_{2}^{\infty}\frac{1}{p^{2}}-\frac{1}{3}\cdot\sum_{2}^{\infty}\frac{1}{p^{3}}-\cdots$$

So \textsc{Mertens} \emph{defines}
\begin{align}
\fbox{$\dis H:=\frac{1}{2}\cdot\sum_{2}^{\infty}\frac{1}{p^{2}}+\frac{1}{3}\cdot\sum_{2}^{\infty}\frac{1}{p^{3}}+\cdots$}
\end{align}Combining this result with (3.3.3) we obtain \begin{align}
\sum_{2}^{\infty}\frac{1}{p^{1+\rho}}=\ln \left(\frac{1}{\rho}\right)-H+o(\rho).
\end{align}which is the equation (3.2) cited earlier.
\subsection{\textsc{Mertens}' Use of Partial Summation}
\textsc{Mertens} wants to compute the \emph{remainder}:$$\sum_{p>x }\frac{1}{p^{1+\rho}}.$$
His object is to show that the ``remainder" series is, effectively, the series $$\sum_{n=G+1}^{\infty}\frac{1}{n^{1+\rho}\ln n},$$ where $G:=[x]$.  That way he reduces his problem to the study of an infinite series over all the integers, something hopefully more amenable to analysis.  He does this by using partial summation.  The form of the partial summation formula which he uses is
\begin{align}
\sum_{n=G+1}^{\infty}a(n)f(n)=\sum_{n=G+1}^{\infty}[A(n)-A(n-1)]f(n)
\end{align}
where he puts:  

\[a(n):= \left\{ \begin{array}{r@{\quad if \quad}l}
\dfrac{\ln p}{p} & n=p\\0 & n\neq n
\end{array}\right.\] and $$f(n):=\frac{1}{n^{\rho}\ln n}.$$

Then, if, with \textsc{Mertens}, we put $G:=[x],$ we perform an almost dizzying sequence of series transformations to obtain:

\begin{align*}
\sum_{p\geqslant G+1}\frac{1}{p^{1+\rho}} & = \sum_{n=G+1}^{\infty}\frac{[A(n)-A(n-1)]}{n^{\rho}\ln n}  \\
 & \overset{(2.1.1)}{=}  -\frac{A(G)}{(G+1)\ln (G+1)}+\sum_{n=G+1}^{\infty}A(n)\left\{\frac{1}{n^{\rho}\ln n}-\frac{1}{(n+1)^{\rho}\ln (n+1)}\right\} \\
&= -\frac{A(G}{(G+1)^{\rho}\ln (G+1)}+
\sum_{n=G+1}^{\infty}\overbrace{\{\ln n+R(n)\}}^{\textstyle{Grossehilfsatz \ 1}}
\left\{\frac{1}{n^{\rho}\ln n}-\frac{1}{(n+1)^{\rho}\ln (n+1)}\right\}\\
&= -\frac{A(G)}{(G+1)^{\rho}\ln (G+1)}+\sum_{n=G+1}^{\infty}R(n)\left\{\frac{1}{n^{\rho}\ln n}-\frac{1}{(n+1)^{\rho}\ln (n+1)}\right\} \\
&\ \  \ \ 
+\sum_{n=G+1}^{\infty}
\left\{\frac{1}{n^{\rho}}-\frac{1}{(n+1)^{\rho}}\underbrace{-\frac{\ln \left(1-\frac{1}{n+1}\right)}{(n+1)^{\rho}\ln (n+1)}}_{\textstyle{=\frac{1}{(n+1)^{1+\rho}\ln (n+1)}+\frac{\lambda}{2n(n+1)^{1+\rho}\ln (n+1)} \ |\lambda|<1}}\right\}\\
&= -\frac{A(G)}{(G+1)^{\rho}\ln (G+1)}+\sum_{n=G+1}^{\infty}R(n)\left\{\frac{1}{n^{\rho}\ln n}-\frac{1}{(n+1)^{\rho}\ln (n+1)}\right\}\\
& \ \ \ \ 
+\sum_{n=G+1}^{\infty}
\left\{\frac{1}{n^{\rho}}-\frac{1}{(n+1)^{\rho}}+\frac{1}{(n+1)^{1+\rho}\ln (n+1)}+\frac{\lambda}{2n(n+1)^{1+\rho}\ln (n+1)} \right\}\\
&= \sum_{n=G+1}^{\infty}\frac{1}{n^{1+\rho}\ln n}+\frac{\ln(G+1)-A(G)}{(G+1)^{\rho}\ln (G+1)}-\frac{1}{(G+1)^{1+\rho}\ln (G+1)}+\\
& \ \ \ \ \lambda\cdot\sum_{n=G+1}^{\infty}\frac{1}{2n(n+1)^{1+\rho}\ln (n+1)}+\\
& \ \ \  +\sum_{n=G+1}^{\infty}R(n)\left\{\frac{1}{n^{\rho}\ln n}-\frac{1}{(n+1)^{\rho}\ln (n+1)}\right\}
\end{align*}and we have proved:
\begin{thm}
\begin{align}
\sum_{p\geqslant G+1}\frac{1}{p^{1+\rho}}=\sum_{n=G+1}^{\infty}\frac{1}{n^{1+\rho}\ln n}+\Re
\end{align}
where
\begin{align*}
\Re&:=\frac{\ln(G+1)-A(G)}{(G+1)^{\rho}\ln (G+1)}-\frac{1}{(G+1)^{1+\rho}\ln (G+1)}+\\
& \ \ \ \ \   \lambda\cdot\sum_{n=G+1}^{\infty}\frac{1}{2n(n+1)^{1+\rho}\ln (n+1)}+\sum_{n=G+1}^{\infty}R(n)\left\{\frac{1}{n^{\rho}\ln n}-\frac{1}{(n+1)^{\rho}\ln (n+1)}\right\}
\end{align*}\hfill$\Box$
\end{thm}
Concerning this rather formidable error term, $\Re$, \textsc{Mertens} writes \emph{``F\"ur $\Re$ es leicht eine obere Grenze anzugeben. . ."} (``It is easy to obtain an upper bound for $\Re$. . .")  He goes on to say that the reason is that by the Grossehilfsatz 1, the numerical value of $R(n)$ can never exceed $2.$  Indeed, as $\rho\rightarrow 0^{+}$ :

\begin{align*}
\frac{\ln(G+1)- A(G)}{(G+1)^{\rho}\ln (G+1)}-\frac{1}{(G+1)^{1+\rho}\ln (G+1)}& =  -\frac{R(G)}{(G+1)^{\rho}\ln \left(G+1\right)}+\\
& \ \ \ \ \ +\frac{\overbrace{\ln \left(1+\dfrac{1}{G}\right)-\dfrac{1}{G+1}}^{\textstyle{<\frac{1}{G^{2}}}}}{(G+1)^{\rho}\ln (G+1)} \\
 & <  \frac{2}{\ln (G+1)}+\frac{1}{
 G^{2}\ln (G+1)}, 
\end{align*}

\noindent and

\begin{align*}
\sum_{n=G+1}^{\infty}\frac{1}{2n(n+1)^{1+\rho}\ln (n+1)} & <  \frac{1}{2}\sum_{n=G+1}^{\infty}\left\{\frac{1}{n\ln n}-\frac{1}{(n+1)\ln (n+1)}\right\} \\
 & =  \frac{1}{2(G+1)\ln (G+1)} 
\end{align*}

\noindent and

\begin{align*}
\sum_{n=G+1}^{\infty}R(n)\left\{\frac{1}{n^{\rho}\ln n}-\frac{1}{(n+1)^{\rho}\ln (n+1)}\right\} & <  2\sum_{n=G+1}^{\infty}\left\{\frac{1}{\ln n}-\frac{1}{\ln (n+1)}\right\} \\
 & =  \frac{2}{\ln (G+1)},
\end{align*}where we used telescopic summation in the last two estimates.  Finally, if $G>2$, then
\begin{align*}
\frac{1}{\ln (G+1)}\left(\frac{1}{G^{2}}+\frac{1}{2(G+1)}\right)&<\frac{1}{\ln (G+1)}\left(\frac{1}{G^{2}}+\frac{1}{2G}\right)\\
&<\frac{1}{\ln (G+1)}\left(\frac{1}{2G}+\frac{1}{2G}\right)\\
&=\frac{1}{G\ln (G+1)}.
\end{align*}

Therefore, we have proved the following error estimate:
\begin{thm}
\begin{align}
|\Re|<\frac{4}{\ln (G+1)}+\frac{1}{G\ln (G+1)}.
\end{align}
\end{thm}\hfill$\Box$
\subsection{Proof the the Grossehilfsatz 1}

We have used the \textbf{Grossehilfsatz 1} on several occasions and the time has come to prove it.  Starting with the standard definition:\begin{align}
\theta(x):=\sum_{p\leqslant x}\ln p,
\end{align}we will use \textsc{Chebyshev}'s technique to prove:
\begin{thm}
\begin{align}\theta(x)<2x.\end{align}
\end{thm}
\begin{proof}
The proof is based on the equation
\begin{align*}
\ln(1\cdot2\cdot3\cdots[x]) & =  \theta(x)+\theta(\sqrt{x})+\theta(\sqrt[3]{x})+\cdots \\
 &   +\theta\left(\frac{x}{2}\right)\theta\left(\sqrt{\frac{x}{2}}\right)\theta\left(\sqrt[3]{\frac{x}{2}}\right)\cdots\\ 
 &  +\theta\left(\frac{x}{3}\right)\theta\left(\sqrt{\frac{x}{3}}\right)\theta\left(\sqrt[3]{\frac{x}{3}}\right)\cdots
 \end{align*}\begin{equation} +\cdots  \ \ \ \ \ \ \ \ \ \ \ \ \ \ \end{equation}

To see why this latter equation is true, define:
\begin{align}
\chi(x):= \theta(x)+\theta(\sqrt{x})+\theta(\sqrt[3]{x})+\cdots.
\end{align} Then we use a well-known theorem of \textsc{Legendre} \cite{Hardy}: \emph{the prime number $p$ divides the number $n!$ exactly 
$$\left[\frac{n}{p}\right]+\left[\frac{n}{p^{2}}\right]+\left[\frac{n}{p^{3}}\right]+\cdots
$$
times.
}Therefore,
\begin{align*}
\ln ([x]!)=\sum_{p\leqslant x}\left(\left[ \frac{x}{p}\right]+\left[ \frac{x}{p^{2}}\right]+\cdots \right)\ln p
\end{align*}

Here, the second member represents the sum of the values of the function $\ln p$ taken over the lattice points $(p,x,u)$, where p is prime, in the region $p>0, \ s>0, \ 0<u\leqslant \frac{x}{p^{s}}$.  The part of the sum which corresponds to two given values of $s$ and $u$ is equal to $\theta\left(\sqrt[x]{\frac{x}{u}}\right)$; the part that corresponds to a given value of $u$ is equal to $\chi\left(\frac{x}{u}\right)$.  

Therefore,
$$\ln(1\cdot2\cdot3\cdots[x])-2\ln\left(1\cdot2\cdot3\cdots\left[\frac{x}{2}\right]\right)=\chi(x)-\chi\left(\frac{x}{2}\right)+\chi\left(\frac{x}{3}\right)-\chi\left(\frac{x}{4}\right)+\cdots.
$$
But, 
$$\chi\left(\frac{x}{3}\right)\geqslant \chi\left(\frac{x}{4}\right), \ \chi\left(\frac{x}{5}\right)\geqslant \chi\left(\frac{x}{6}\right), \ \cdots
$$
and therefore
$$\chi(x)-\chi\left(\frac{x}{2}\right)<\ln(1\cdot2\cdot3\cdots[x])-2\ln\left(1\cdot2\cdot3\cdots\left[\frac{x}{2}\right]\right).
$$
Applying Stirling's formula (3.2.2) and (3.2.3) we obtain that for all $x\geqslant 4$:
\begin{align*}
\chi(x)-\chi\left(\frac{x}{2}\right)&<x\ln 2+\frac{3}{2}\ln x-\ln \sqrt{2\pi}-\ln 2+\frac{2}{x-2}+\frac{1}{12x}\\
&< x-\left\{(1-\ln 2)x-\frac{3}{2}\ln x+\ln\sqrt{2\pi}+\ln 2-\frac{2}{x-2}-\frac{1}{12x}   \right\}\\
&<x
\end{align*}
But this same inequality can be verified directly for $x<4.$  Therefore, we have proved the general inequality:  \emph{if $x>1$, then}
\begin{align}
\chi(x)-\chi\left(\frac{x}{2}\right)<x.
\end{align}
We now substitute $x,\dis \ \frac{x}{2}, \ \frac{x}{4}, \ \frac{x}{8}, \cdots$ for $x$ until we reach a term $\dfrac{x}{2^{m}}$ which is less than $2$.  We then add up the inequalities
\begin{align*}
 \chi(x)-\chi\left(\frac{x}{2}\right)& <  x \\
\chi\left(\frac{x}{2}\right)-\chi\left(\frac{x}{4}\right) & <  \frac{x}{2}\\ 
\chi\left(\frac{x}{4}\right)-\chi\left(\frac{x}{8}\right) & <  \frac{x}{4}\\
........................... & <  ......
\end{align*}and we obtain
\begin{align*}
\chi(x)&<x\left(1+\frac{1}{2}+\frac{1}{4}+\cdots +\frac{1}{2^{m}}\right)\\
& <  2x,
\end{align*}
and so all the more is
$$\theta(x)<2x$$
\end{proof}

Chebyshev, himself, proved \cite{Chebyshev} that $$0.904x<\theta(x)<1.113x$$ for $x\geqslant 38750.$

Now we are ready to complete the proof of the \textbf{Grossehilfsatz 1}.  We use the inequality for $\theta(x)$ and Legendre's theorem again. This latter implies that
$$\ln n!=\sum_{p\leqslant n}\left[\frac{n}{p}\right]\ln p+\sum_{p^{2}\leqslant n}\left[\frac{n}{p^{2}}\right]\ln p+\sum_{p^{3}\leqslant n}\left[\frac{n}{p^{3}}\right]\ln p+\cdots.
$$
If we write
$$\left[\frac{n}{p}\right]:=\frac{n}{p}-r_{p},
$$
and use Stirling's formula (3.2.4), we obtain
\begin{align}
\ln n-1+\frac{1}{2n}\ln n+\frac{\ln \sqrt{2\pi}}{n}+\frac{\lambda}{12n^{2}}=\sum_{p\leqslant n}\frac{\ln p}{p}-\frac{1}{n}\sum_{p\leqslant n}r_{p}\ln p+\frac{1}{n}\sum_{p^{2}\leqslant n}\left[\frac{n}{p^{2}}\right]\ln p +\cdots.
\end{align}Here, $|\lambda|<1.$  We rewrite this as:
\begin{align}
\ln n-\sum_{p\leqslant n}\frac{\ln p}{p}=1-\frac{1}{2n}\ln n-\frac{\ln \sqrt{2\pi}}{n}-\frac{\lambda}{12n^{2}}-\frac{1}{n}\sum_{p\leqslant n}r_{p}\ln p+\frac{1}{n}\sum_{p^{2}\leqslant n}\left[\frac{n}{p^{2}}\right]\ln p +\cdots.
\end{align}

Therefore, if $n\geqslant 5$, the equation (3.5.7) shows that $\left\{\dis\ln n -\sum_{p\leqslant n}\frac{\ln p}{p}\right\}$ is contained between the \emph{upperbound}
$$1+\sum_{p^{2}\leqslant n}\frac{\ln p}{p^{2}}+\sum_{p^{3}\leqslant n}\frac{\ln p}{p^{3}}+\cdots
$$
and the \emph{lower bound}
$$-\frac{1}{n}\sum_{p\leqslant n}\frac{\ln p}{p}.
$$
Now, on the one hand, by \textbf{Theorem 11},
$$\sum_{p\leqslant n}\ln p<2n,$$ 
while, on the other hand,

\begin{align*}
 \sum_{p^{2}\leqslant n}\frac{\ln p}{p^{2}}+\sum_{p^{3}\leqslant n}\frac{\ln p}{p^{3}}+\cdots
 & <  \sum_{p\geqslant 2}^{\infty}\frac{\ln p}{p^{2}}+\sum_{p\geqslant 2}^{\infty}\frac{\ln p}{p^{3}}+\sum_{p\geqslant 2}^{\infty}\frac{\ln p}{p^{4}}+\sum_{p\geqslant 2}^{\infty}\frac{\ln p}{p^{5}}+\cdots \\
& <  \sum_{p\geqslant 2}^{\infty}\frac{\ln p}{p^{2}}+\frac{1}{2}\sum_{p\geqslant 2}^{\infty}\frac{\ln p}{p^{2}}+\sum_{p\geqslant 2}^{\infty}\frac{\ln p}{p^{4}}+\frac{1}{2}\sum_{p\geqslant 2}^{\infty}\frac{\ln p}{p^{4}}+\cdots \\
& = \frac{3}{2}\left(\sum_{p\geqslant 2}^{\infty}\frac{\ln p}{p^{2}}+\sum_{p\geqslant 2}^{\infty}\frac{\ln p}{p^{4}}+\cdots\right) \\
  & =  \frac{3}{2} \sum_{p\geqslant 2}^{\infty}\frac{\ln p}{p^{2}}\left(1+\frac{1}{p^{2}}+\frac{1}{p^{4}}+\cdots \right)\\
   &=\frac{3}{2}\sum_{p\geqslant 2}^{\infty}\left\{\frac{\dis\frac{\ln p}{p^{2}}}{1-\dfrac{1}{p^{2}}}\right\}\\  
 &= \frac{3}{2}\left\{\frac{\dis\sum_{n=1}^{\infty}\frac{\ln n}{n^{2}}}{\dis\dis\sum_{n=1}^{\infty}\frac{1}{n^{2}}}\right\}\ \ 
\\&=\frac{3}{2}\cdot \frac{0.9375482543...}{\dfrac{\pi^{2}}{6}}<\frac{9}{\pi^{2}}<1.
\end{align*}The penultimate equality is the logarithmic derivative of \textsc{Euler}'s identity at $\rho=1.$
Therefore, we have proven that for $n>4,$
$$\left|\dis\ln n -\sum_{p\leqslant n}\frac{\ln p}{p}\right|<2
$$
Finally, for $1\leqslant n\leqslant4$, (see \cite{Bachmann})
$$1-\frac{\ln \sqrt{2\pi n}}{2}-\frac{\lambda}{12n^{2}}>0
$$
because
$$\frac{\ln \sqrt{2\pi n}}{n}=\frac{\ln 2n}{2n}+\frac{\ln \pi}{2n}<\frac{\ln 2}{4}+\frac{\ln 2}{4}=\ln 2
$$and
$$\frac{\lambda}{12n^{2}}<\frac{1}{48}
$$
and therefore,
$$\frac{\ln \sqrt{2\pi n}}{2}+\frac{\lambda}{12n^{2}}<\ln 2+\frac{1}{48}<1.
$$
This completes the proof of the \textbf{Grossehilfsatz 1}.
\hfill$\Box$

The reader will observe that the more accurate inequality of \textsc{Chebyshev}, $\theta(x)<1.13x$ is of no use in improving the bound which \textsc{Mertens} obtained in the \textbf{Grossehilfsatz 1}, since it is used to obtain the \emph{lower} bound, only, while the \emph{upper} bound is of the form $1+(1-\epsilon)$ where $\epsilon$ is very tiny, and for which the results of \textsc{Chebyshev} are irrelevant.
Using the most advanced techniques available, \textsc{Dusart} \cite{Dus} has proven:
$$\lim_{x\rightarrow \infty}\left\{\sum_{p\leqslant x}\frac{\ln p}{p}-\ln x\right\}=-1.3325822757...
$$
So the value $2$ given by \textsc{Mertens} as an upper bound for the absolute value of the constant is pretty close to the true value.

\subsection{The Grossehilfsatz 2}
We state:

\begin{lemma}
\begin{align}
\fbox{$\dis\sum_{n= G+1}^{\infty}\frac{1}{n^{1+\rho}\ln n}=\ln\ln G+\gamma+ \frac{\lambda}{G\ln G}+o(\rho). $}
\end{align}
where $\gamma$ is \textsc{Euler}'s constant, and $|\lambda|< 1.$
\end{lemma}We offer \emph{two} proofs.  \textsc{Mertens}' original proof, which displays his technical virtuosity, and our own modern proof.


\subsubsection{Merten's proof}
\begin{proof}
This is another marvelous \emph{tour de force}. 

The first step is to obtain an estimate for the ``remainder" in the \textsc{Riemann} zeta-function: $\sum_{n=G+1}^{\infty}\frac{1}{n^{1+\rho}}$.

We begin by noting that the binomial theorem gives us

\begin{align*}
\frac{1}{n^{\rho}} & =  \frac{1}{(n+1)^{\rho}}\left(\frac{n+1}{n}\right)^{\rho} \\
 & =   \frac{1}{(n+1)^{\rho}}\left(\frac{1}{1-\dfrac{1}{n+1}}\right)^{\rho}\\
& =  \frac{1}{(n+1)^{\rho}}\left(1-\frac{1}{n+1}\right)^{-\rho}\\
& =  \frac{1}{(n+1)^{\rho}}\left\{1+\frac{\rho}{1!}\frac{1}{(n+1)}+\frac{\rho(\rho+1)}{2!}\frac{1}{(n+1)^{2}}+\frac{\rho(\rho+1)(\rho+2)}{3!}\frac{1}{(n+1)^{3}}+\cdots\right\}\\
& =  \frac{1}{(n+1)^{\rho}}+\frac{\rho}{1!}\frac{1}{(n+1)^{1+\rho}}+\frac{\rho(\rho+1)}{2!}\frac{1}{(n+1)^{2+\rho}}+\frac{\rho(\rho+1)(\rho+2)}{3!}\frac{1}{(n+1)^{3+\rho}}+\cdots
\end{align*}
and transposing the first term on the left to the right hand side and dividing both sides by $\rho$ we obtain:
$$\frac{1}{\rho n^{\rho}}-\frac{1}{\rho(n+1)^{\rho}}=\frac{1}{1!}\frac{1}{(n+1)^{1+\rho}}+\frac{(\rho+1)}{2!}\frac{1}{(n+1)^{2+\rho}}+\frac{(\rho+1)(\rho+2)}{3!}\frac{1}{(n+1)^{3+\rho}}+\cdots
$$
If we sum this last equation from $n=G$ to $n=\infty$ we obtain:
\begin{align}
\sum_{n=G+1}^{\infty}\frac{1}{n^{1+\rho}}=\frac{1}{\rho G^{\rho}}-\Re'
\end{align}
where
\begin{align}
\Re'=\frac{(\rho+1)}{2!}\sum_{n=G+1}^{\infty}\frac{1}{(n+1)^{2+\rho}}+\frac{(\rho+1)(\rho+2)}{3!}\sum_{n=G+1}^{\infty}\frac{1}{(n+1)^{3+\rho}}+\cdots
\end{align}

We have now obtained the promised representation of the ``remainder." The next step is as marvelous as it is unexpected.  We integrate (3.6.2) with respect to \emph{the exponent}, $\rho$\ !

The summand, $\dis\frac{1}{n^{1+\rho}\ln n}$, can be obtained from the identity:
$$\int_{\rho}^{1}\frac{1}{n^{1+t}}~dt=\frac{1}{n^{1+\rho}\ln n}-\frac{1}{n^{2}\ln n}
$$

If we apply this to (3.6.2) and (3.6.3) by integrating them from $t=\rho$ to $t=1$ we obtain
\begin{align*}
&\ \ \ \ \ \ \ \ \ \ \ \ \ \ \ \ \ \ \ \ \ \ \ \ \sum_{n= G+1}^{\infty}\frac{1}{n^{1+\rho}\ln n}-\sum_{n= G+1}^{\infty}\frac{1}{n^{2}\ln n}=\\
&=\int_{\rho}^{1}\frac{1}{tG^{t}}~dt-\int_{\rho}^{1}\Re'~dt \\
& =  \int_{\rho}^{\infty}\frac{1}{tG^{t}}~dt-\int_{1}^{\infty}\frac{1}{tG^{t}}~dt-\int_{\rho}^{1}\Re'~dt\\
& = \underbrace{ \int_{\rho \ln G}^{\infty}\frac{1}{xe^{x}}~dx}_{\textstyle{x:=t\ln G}}-\int_{1}^{\infty}\frac{1}{tG^{t}}~dt-\int_{\rho}^{1}\Re'~dt\\
& = \underbrace{ \int_{\rho \ln G}^{\infty}\frac{1}{e^{x}-1}~dx}_{\textstyle{=-\ln(1-\frac{1}{G^{\rho}}})}-\int_{\rho \ln G}^{\infty}\left\{\frac{1}{e^{x}-1}-\frac{1}{xe^{x}}\right\}~dx
  -\int_{1}^{\infty}\frac{1}{tG^{t}}~dt-\int_{\rho}^{1}\Re'~dt\\ 
& = -\ln\left(1-\frac{1}{G^{\rho}}\right)-\underbrace{\int_{0}^{\infty}\left\{\frac{1}{e^{x}-1}-\frac{1}{xe^{x}}\right\}~dx}_{\textstyle{= \ \gamma \  (\rm \textsc{Euler}'s\  \rm constant)}}
  +\underbrace{\int_{0}^{\rho\ln G}\left\{\frac{1}{e^{x}-1}-\frac{1}{xe^{x}}\right\}~dx}_{\textstyle{< \ \rho\ln G \ \rm if\  \rho<\ln G^{2}}}\\ 
&\ \ \ \ \ -\int_{1}^{\infty}\frac{1}{tG^{t}}~dt-\int_{\rho}^{1}\Re'~dt\\
& =  \ln \left(\frac{1}{\rho}\right)-\ln\ln G-\gamma -\int_{1}^{\infty}\frac{1}{tG^{t}}~dt-\int_{\rho}^{1}\Re'~dt+o(\rho),\\
\end{align*}since
\begin{align*}
 -\ln\left(1-\frac{1}{G^{\rho}}\right)&=-\ln(1-e^{-\rho\ln G})\\
 &=-\ln(1-\{1-\rho\ln G+o(\rho)\})\\
 &=-\ln \rho -\ln\ln G+o(\rho)\\
 &=\ln\left(\frac{1}{\rho}\right)-\ln\ln G +o(\rho).
\end{align*}
and therefore,
\begin{align}
\sum_{n= G+1}^{\infty}\frac{1}{n^{1+\rho}\ln n}=\ln \left(\frac{1}{\rho}\right)-\ln\ln G-\gamma \underbrace{-\int_{1}^{\infty}\frac{1}{tG^{t}}~dt+\sum_{n= G+1}^{\infty}\frac{1}{n^{2}\ln n}-\int_{\rho}^{1}\Re'~dt}_{\textstyle{=\epsilon \  \equiv\  \rm error}}+o(\rho)
\end{align}

This shows where the \textsc{Euler}'s constant component of \textsc{Mertens}' constant $B$ comes from.  Namely, from a subtle and delicate trick of \emph{adding and subtracting the nonobvious integral $\int_{\rho \ln G}^{\infty}\frac{1}{e^{x}-1}~dx$ to and from the sum $\sum_{n= G+1}^{\infty}\frac{1}{n^{1+\rho}\ln n}$.}

Now we estimate the error:

\begin{align*}
\int_{\rho}^{1}\Re'~dt & <  \int_{0}^{1}\left\{\sum_{n= G+1}^{\infty}\frac{1}{n^{2+t}}+\sum_{n= G+1}^{\infty}\frac{1}{n^{3+t}}+\cdots \right\}~dt\\
 & <   \sum_{n= G+1}^{\infty}\left(\frac{1}{n^{2}\ln n}-\frac{1}{n^{3}\ln n}\right)+\sum_{n= G+1}^{\infty}\left(\frac{1}{n^{3}\ln n}-\frac{1}{n^{4}\ln n}\right)+\cdots\\
& =  \sum_{n= G+1}^{\infty}\frac{1}{n^{2}\ln n}\\
& <  \sum_{n= G+1}^{\infty}\left\{\frac{1}{(n-1)\ln (n-1)}-\frac{1}{n\ln n}\right\}\\
& =  \frac{1}{G\ln G},
\end{align*}
and

$$\int_{1}^{\infty}\frac{1}{tG^{t}}~dt<\int_{1}^{\infty}\frac{1}{G^{t}}~dt=\frac{1}{G\ln G}.
$$
Therefore,

\begin{align*}
\epsilon & =  -\int_{1}^{\infty}\frac{1}{tG^{t}}~dt+\sum_{n= G+1}^{\infty}\frac{1}{n^{2}\ln n}-\int_{\rho}^{1}\Re'~dt\\
& =  \lambda_{1}\sum_{n= G+1}^{\infty}\frac{1}{n^{2}\ln n}-\frac{\lambda_{2}}{G\ln G}\\
& =  \frac{\lambda_{3}-\lambda_{2}}{G\ln G}\\
& <  \frac{1}{G\ln G}
\end{align*}
where $0<\lambda_{k}<1$ for $k=1, \ 2, \ 3.$  

We have shown:

$$\sum_{n= G+1}^{\infty}\frac{1}{n^{1+\rho}\ln n}=\ln\left(\frac{1}{\rho}\right)-\ln\ln G-\gamma +\frac{\lambda}{G\ln G}+o(\rho).
$$
where $|\lambda|<1.$ 
This completes the proof of the \textbf{Grossehilfsatz 2}. \hfill$\Box$

\subsubsection{Modern Proof }
It may be of interest to insert a modern proof of \textbf{Grossehilfsatz 2} based on a simple form of the \textsc{Euler-MacLaurin} formula as given by \textsc{Boas\ \cite{Boas} .}
\begin{thm}
Let $f(t)$ be positive for $t>0$ and suppose that $|f'(t)|$ is decreasing.  If $\sum_{n=1}^{\infty}f(n)$ is convergent and if $$R_{n}:=f(n+1)+f(n+2)+\cdots,$$ then there exists a number $\theta$ with $0<\theta<1$ such that the following equation is valid:
\begin{align}
R_{n}=\int_{n+\frac{1}{2}}^{\infty}f(t)~dt+\frac{\theta}{8}f'(n+1).
\end{align}

\end{thm}\hfill$\Box$

In the coming computation, we will use the following results.
For fixed $G$,
\begin{align}
\left(G+\frac{1}{2}\right)^{-\rho}=(G+1)^{-\rho}=1+o(\rho)
\end{align}
since, for any contant, $\alpha$, $$(G+\alpha)^{-\rho}=e^{-\rho\ln(G+\alpha)}=1+\rho\ln(G+\alpha)-\frac{1}{2}\{\rho\ln(G+\alpha)\}^{2}+\cdots=1+o(\rho)$$
Moreover, by \textsc{Taylor}'s theorem 
\begin{align}
\ln(1+x)=x-\frac{\lambda}{2}x^{2}.
\end{align}
where $0<\lambda<1.$
Finally,
\begin{align}
-\gamma=\int_{0}^{\infty}\frac{\ln v}{e^{v}}~dv
\end{align}
which follows from the change of variable $x:=e^{v}$ in the standard integral $$-\gamma=\int_{0}^{1}\ln\ln\frac{1}{x}~dx,$$
which appears in \textsc{Havil \ \cite{Havil}}, p. 109.

Then, substituting in (3.6.5) and integrating by parts with $$u:= x^{-\rho}, \ dv:=\frac{dx}{x\ln x},$$we obtain

\begin{align*} 
&\ \ \ \ \ \ \ \ \ \ \ \ \ \ \ \ \ \ \ \ \ \ \ \ \ \sum_{n= G+1}^{\infty}\frac{1}{n^{1+\rho}\ln n}  =\\ 
&=\int_{G+\frac{1}{2}}^{\infty}\frac{1}{x^{1+\rho}\ln x}~dx+\frac{\theta}{8}\left\{\frac{1}{x^{1+\rho}\ln x}\right\}'_{x=G+1}  \\
 &=  \frac{\ln\ln x}{x^{\rho}}\Bigg|_{x=G+\frac{1}{2}}^{x=\infty}-\int_{G+\frac{1}{2}}^\infty\frac{(\ln\ln x)(-\rho)}{x^{\rho +1}}~dx
 - \frac{\theta}{8(G+1)^{2+\rho}}\left\{1+\rho +\frac{1}{\ln(G+1)}\right\} \\
&=- \frac{\ln\ln(G+\frac{1}{2})}{(G+1)^{\rho}}+\rho\int_{G+\frac{1}{2}}^\infty\frac{(\ln\ln x)}{x^{\rho +1}}~dx
 - \frac{\theta}{8(G+1)^{2+\rho}}\left\{1+\rho +\frac{1}{\ln(G+1)}\right\} \\
&\overset{(3.6.6)}{=}-\ln\ln\left(G+\frac{1}{2}\right)+\rho\int_{G+\frac{1}{2}}^\infty\frac{(\ln\ln x)}{x^{\rho +1}}~dx
 - \frac{\theta}{8(G+1)^{2}}\left\{1 +\frac{1}{\ln(G+1)}\right\}+o(\rho)\\
&=-\ln\ln\left(G+\frac{1}{2}\right)+\rho\int_{G+\frac{1}{2}}^\infty\frac{(\ln\ln x)}{x^{\rho +1}}~dx
 - \frac{\theta_{1}}{4(G+1)^{2}}+o(\rho)\ \ \ \ \ (0<\theta_{1}<1)\\
&=-\ln\ln G-\ln\left\{1+\frac{\ln\left(1+\frac{1}{2G}\right)}{\ln G}\right\}+\rho\int_{G+\frac{1}{2}}^\infty\frac{(\ln\ln x)}{x^{\rho +1}}~dx
 - \frac{\theta_{1}}{4(G+1)^{2}}+o(\rho)\\
&\overset{(3.6.7)}{=}-\ln\ln G-\frac{\theta_{2}}{2G\ln G}+\rho\int_{G+\frac{1}{2}}^\infty\frac{(\ln\ln x)}{x^{\rho +1}}~dx
 - \frac{\theta_{1}}{4(G+1)^{2}}+o(\rho)\ \ (0<\theta_{2}<1)\\
&=-\ln\ln G+\rho\int_{G+\frac{1}{2}}^\infty\frac{(\ln\ln x)}{x^{\rho +1}}~dx
 -\frac{\theta_{3}}{G\ln G}+o(\rho)\ \ \ \ (0<\theta_{3}<1)\\
&\overset{(x:=e^{\frac{v}{\rho}})}{=}-\ln\ln G+\int_{\rho\ln(G+\frac{1}{2})}^{\infty}\frac{\ln\frac{1}{\rho}}{e^{v}}~dv+\int_{\rho\ln(G+\frac{1}{2})}^{\infty}\frac{\ln v}{e^{v}}~dv
 -\frac{\theta_{3}}{G\ln G}+o(\rho)\\
&=-\ln\ln G+\frac{\ln \frac{1}{\rho}}{(G+\frac{1}{2})^{\rho}}+\int_{\rho\ln(G+\frac{1}{2})}^{\infty}\frac{\ln v}{e^{v}}~dv-\frac{\theta_{3}}{G\ln G}+o(\rho)\\
&\overset{(3.6.6)}{=}\ln\frac{1}{\rho}-\ln\ln G+\int_{\rho\ln(G+\frac{1}{2})}^{\infty}\frac{\ln v}{e^{v}}~dv-\frac{\theta_{3}}{G\ln G}+o(\rho)\\
&=\ln\frac{1}{\rho}-\ln\ln G+\int_{0}^{\infty}\frac{\ln v}{e^{v}}~dv-\int_{0}^{\rho\ln(G+\frac{1}{2})}\frac{\ln v}{e^{v}}~dv-\frac{\theta_{3}}{G\ln G}+o(\rho)\\
&\overset{(3.6.8)}{=}\ln\frac{1}{\rho}-\ln\ln G-\gamma+o(\rho)-\frac{\theta_{3}}{G\ln G}+o(\rho)\\
&=\ln\frac{1}{\rho}-\ln\ln G-\gamma-\frac{\theta_{3}}{G\ln G}+o(\rho)\\
\end{align*}\hfill$\Box$

Observe that this method produces the dominant terms $$\ln\frac{1}{\rho}, \ -\ln\ln G,\ \rm \textsc{Euler}'s\  constant=\gamma,$$ almost automatically, without the nonobvious and tricky (but beautiful and clever) artifices employed by \textsc{Mertens},  while the error term, $-\frac{\theta_{3}}{G\ln G}$, with a \emph{sign}, appears with virtually no effort.  The reason is the power of the half-interval version of the \textsc{Euler-Maclaurin} formula combined with the use of integration by parts.  I think that \textsc{Mertens} would have liked this proof.

\subsection{The Formula for the Constant $H$}
\textsc{Mertens} computes the constant $B:=\gamma-H$ by finding a rapidly convergent series for $H.$  The paper \cite{L-P} treats the computation exhaustively. However, they do not give \textsc{Mertens}' own derivation, so we develop it here.  Define:
$$x_{k}:=\frac{1}{k}\sum_{p\geqslant 2}^{\infty}\frac{1}{p^{k}}, \ \zeta(k):=\sum_{n=1}^{\infty}\frac{1}{n^{k}}.
$$

Then, by (3.3.4)

\begin{align}
H  &=&  x_{2}+x_{3}+x_{4}+x_{5}+x_{6}+x_{7}+x_{8}+\cdots\\
\frac{1}{2}\ln\{\zeta(2)\} & = & x_{2}+ \ \ \ \ + x_{4}+ \ \ \ \ + x_{6}+ \ \ \ \ +\ x_{8}+\cdots\\
\frac{1}{3}\ln\{\zeta(3)\} & = &\   \ \ \ \ \ \ x_{3}+ \ \ \ \ \ \ \ \  \ \ + x_{6}+ \ \ \ \ \ \ \ \ \ \ \  +\cdots\\
\frac{1}{4}\ln\{\zeta(4)\} & = & \ \ \ \ \ \  \ \ \ \ \  +x_{4}\ \  \ \ \ \ \ \ \ \ \ \ \ \ \ \ \ \  \ \ + x_{8}+\cdots
\end{align}

and so on.  Now, let $\mu(n)$:
\begin{enumerate}
  \item have the value $1$, if $n=1$, or has an even number of distinct prime divisors.
  \item have the value $-1$ if $n$ has an odd number of distinct prime divisors.
  \item vanish if $n$ is equal to a prime divisor.
\end{enumerate}
 
 Moreover, let $1,\  d,\  d',\  \cdots $ be all the divisors of $n$.  Then it follows from the definition of the numbers $\mu(1), \ \mu(2), \ \mu(3), \ \cdots$, that for any integer $n$ greater than $1$, 
 \begin{align}
 \mu(1)+\mu(d)+\mu(d')+\cdots =0
 \end{align}

Now, if we multiply the equations (3.7.1), (3.7.2), (3.7.3), etc. by $\mu(1)$, $\mu(2)$, $\mu(3)$, etc., respectively, and add up the resulting equations and use (3.7.5), we see that $x_{1}$, $x_{2}$, $x_{3}$, ... all drop out and we obtain: 
$$H-\frac{1}{2}\ln\{\zeta(2)\}-\frac{1}{3}\ln\{\zeta(3)\}-\frac{1}{5}\ln\{\zeta(5)\}+\frac{1}{6}\ln\{\zeta(6)\}-\frac{1}{7}\ln\{\zeta(7)\}+\frac{1}{10}\ln\{\zeta(10)\}-\cdots=0
$$

Therefore, he have proved:
\begin{thm} 
 \begin{align*}
\fbox{$\dis H=\frac{1}{2}\ln\{\zeta(2)\}+\frac{1}{3}\ln\{\zeta(3)\}+\frac{1}{5}\ln\{\zeta(5)\}-\frac{1}{6}\ln\{\zeta(6)\}+\frac{1}{7}\ln\{\zeta(7)\}-\frac{1}{10}\ln\{\zeta(10)\}+\cdots$}
\end{align*}
\end{thm}\hfill$\Box$

We observe that the absolute convergence of the series in question allow the elimination of the $x_{k}$'s.  Using the published tables of \textsc{Legendre} \cite{L} of the values of $\zeta(m)$ to fifteen decimal places, \textsc{Mertens} computed the value:

$$H\approx 0.31571845205,$$

\noindent and therefore,

$$B=\gamma-H\approx 0.2614972128.
$$


\subsection{Completion of the Proof}

Now we follow the sketch in 3.1.  

\begin{align*}
 \sum_{p\leqslant x}\frac{1}{p^{1+\rho}}&=\sum_{p}\frac{1}{p^{1+\rho}}-\sum_{p\geqslant x}\frac{1}{p^{1+\rho}} \\
 & =  \ln\left(\frac{1}{\rho}\right)-H+o(\rho)-\sum_{p\geqslant x}\frac{1}{p^{1+\rho}} \ \ \ \ \ \ \ \ (\rm by\  (3.3.5))  \\
& =  \ln\left(\frac{1}{\rho}\right)-H+o(\rho)-\sum_{n=G+1}^{\infty}\frac{1}{n^{1+\rho}\ln n}-\Re \ \ \ \ \ \ \ \ (\rm by\  (3.4.2))\\
& =  \ln\ln G+\gamma-H+\frac{\lambda}{G\ln G}-\Re +o(\rho) \ \ \ \ \ \ \ \ \hfill(\rm by\  (3.6.1))\\
& =  \ln\ln G+\gamma-H+\delta +o(\rho),\ \ \ \ \ \ \ \ (\rm by (3.4.3))
\end{align*}where
$$\fbox{$\dis|\delta|<\frac{4}{\ln(G+1)}+\frac{2}{G\ln G}. $}
$$

Letting $\rho \rightarrow 0$ we obtain
$$\fbox{$\dis\sum_{p\leqslant x}\frac{1}{p^{1+\rho}}=\dis \ln\ln G+\gamma-H+\delta.$}
$$
This completes \textsc{Mertens}' proof of \textsc{Mertens}' Theorem.

\end{proof}
\section{Retrospect and Prospect}
\subsection{Retrospect}

Is this proof not stunning?  The basic idea, totally different from the modern method, is to work with the convergent ``prime zeta function" and study the remainder as $\rho\rightarrow 0^{+}$.  The modern proof is a direct use of partial summation on the given sum.

\textsc{Mertens}' proof is quite natural in approach, and the constant $H$ appears quite inevitably.  The series computations and the manipulation of inequalities are breathtaking.  His use of partial summation is brilliant; indeed, it was hailed as a new technique in prime number theory by contemporaries \cite{Bachmann}.  Finally we signal the repeated clever use of telescopic summations in the estimation of error terms.

Any contemporary analyst can marvel at and be instructed by \textsc{Mertens}' ``arabesques of algebra," a telling phrase due to \textsc{E.T. Bell} \cite{Bell} to describe the manipulations of \textsc{Jacobi} in the theory of elliptic functions to discover number-theoretic theorems, but equally applicable to \textsc{Mertens}' mathematics in this memoir.  

All the techniques \textsc{Mertens} used are now standard tools for the analytic number theorist (among others), but it is a joy to see them used together in a single focused effort to obtain his one towering result.

\subsection{Prospect}

Modern work on \textsc{Mertens}' theorem has concentrated on improving the error term.  The best result to date which has been completely proven is due to \textsc{Dusart} \cite{Dus}:
\begin{thm} For $x>1$
\begin{align*}
\sum_{p\leqslant x}\frac{1}{p}-\ln\ln x -B\geqslant -\left(\frac{1}{10\ln^{2}x}+\frac{4}{15\ln^{3}x}\right)
\end{align*} For $x\geqslant 10372$

\begin{align*}
\sum_{p\leqslant x}\frac{1}{p}-\ln\ln x -B\leqslant \left(\frac{1}{10\ln^{2}x}+\frac{4}{15\ln^{3}x}\right)
\end{align*}\hfill$\Box$

\end{thm}
The best result to date,\emph{ assuming the validity of the \textsc{Riemann} Hypothesis (!)}, is due to \textsc{Schoenfeld} \cite{Schoenfeld}, and affirms:

\begin{thm}
If $x\geqslant 13.5$, then:
\begin{align*}
\fbox{$\left|\dis \sum_{p\leqslant x}\frac{1}{p}-\ln\ln x-B\right|<\dfrac{3\ln x+4}{8\pi\sqrt{x}}$}
\end{align*}
\end{thm}\hfill$\Box$

In both cases, the error term is much better than that of \textsc{Mertens}, himself, but no optimal error term has been found.

Recently, \textsc{M. Wolf} \cite{Wolf1} derived \textsc{Mertens}' series by a completely different method.  He uses the ``generalized \textsc{Bruns} constants" which measure the gaps between consecutive primes, and by an ingenious combination of hard rigorous computations and heuristic numerical arguments obtains \textsc{Mertens}' series, including the big ``O" error term.  Moreover, he prepared a numerical table (which I reproduce with his permission) comparing the error term in \textbf{Theorem 15} with the true error.
\begin{center}
{\sf The Ratio of the True Error to the Predicted Error}\\
\bigskip  
\
\begin{tabular}{|c|c|c|c|} \hline
$x$ & $ |\sum_{p<x} 1/p-\log\log x -B | $ & $\frac{3\log(x)+4}{8\pi\sqrt{x}}$  &  ratio column 2/ column 3 \\ \hline
$ 2^{16}$=          65536&  2.43328226E-0004 &  5.79284588E-0003 &  4.20049542E-0002 \\ \hline
$ 2^{17}$=         131072&  2.26479291E-0004 &  4.32469516E-0003 &  5.23688450E-0002 \\ \hline
$ 2^{18}$=         262144&  1.11367788E-0004 &  3.21961962E-0003 &  3.45903559E-0002 \\ \hline
$ 2^{19}$=         524288&  1.23916030E-0004 &  2.39088215E-0003 &  5.18285814E-0002 \\ \hline
$ 2^{20}$=        1048576&  5.58449145E-0005 &  1.77140815E-0003 &  3.15257184E-0002 \\ \hline
$ 2^{21}$=        2097152&  4.63383665E-0005 &  1.30970835E-0003 &  3.53806756E-0002 \\ \hline
$ 2^{22}$=        4194304&  3.20736392E-0005 &  9.66503244E-0004 &  3.31852370E-0002 \\ \hline
$ 2^{23}$=        8388608&  1.83353157E-0005 &  7.11987819E-0004 &  2.57522885E-0002 \\ \hline
$ 2^{24}$=       16777216&  1.10324946E-0005 &  5.23651207E-0004 &  2.10684030E-0002 \\ \hline
$ 2^{25}$=       33554432&  1.29876787E-0005 &  3.84560730E-0004 &  3.37727640E-0002 \\ \hline
$ 2^{26}$=       67108864&  6.42047777E-0006 &  2.82025396E-0004 &  2.27656015E-0002 \\ \hline
$ 2^{27}$=      134217728&  3.69019851E-0006 &  2.06563775E-0004 &  1.78646934E-0002 \\ \hline
$ 2^{28}$=      268435456&  3.19579180E-0006 &  1.51112594E-0004 &  2.11484146E-0002 \\ \hline
$ 2^{29}$=      536870912&  1.63321145E-0006 &  1.10423592E-0004 &  1.47904212E-0002 \\ \hline
$ 2^{30}$=     1073741824&  1.72440466E-0006 &  8.06062453E-0005 &  2.13929411E-0002 \\ \hline
$ 2^{31}$=     2147483648&  8.53875133E-0007 &  5.87826489E-0005 &  1.45259723E-0002 \\ \hline
$ 2^{32}$=     4294967296&  5.34863712E-0007 &  4.28280967E-0005 &  1.24886173E-0002 \\ \hline
$ 2^{33}$=     8589934592&  5.56640268E-0007 &  3.11767507E-0005 &  1.78543387E-0002 \\ \hline
$ 2^{34}$=    17179869184&  3.62687244E-0007 &  2.26765354E-0005 &  1.59939443E-0002 \\ \hline
$ 2^{35}$=    34359738368&  1.45653226E-0007 &  1.64810885E-0005 &  8.83759748E-0003 \\ \hline
$ 2^{36}$=    68719476736&  1.16826187E-0007 &  1.19695112E-0005 &  9.76031397E-0003 \\ \hline
$ 2^{37}$=   137438953472&  9.94572329E-0008 &  8.68690083E-0006 &  1.14491042E-0002 \\ \hline
$ 2^{38}$=   274877906944&  6.52601557E-0008 &  6.30037736E-0006 &  1.03581344E-0002 \\ \hline
$ 2^{39}$=   549755813888&  5.50727125E-0008 &  4.56662870E-0006 &  1.20598183E-0002 \\ \hline
$ 2^{40}$=  1099511627776&  3.20547296E-0008 &  3.30799956E-0006 &  9.69006466E-0003 \\ \hline
$ 2^{41}$=  2199023255552&  1.70901151E-0008 &  2.39490349E-0006 &  7.13603497E-0003 \\ \hline
$ 2^{42}$=  4398046511104&  1.95113765E-0008 &  1.73290522E-0006 &  1.12593442E-0002 \\ \hline
$ 2^{43}$=  8796093022208&  9.40614690E-0009 &  1.25324631E-0006 &  7.50542552E-0003 \\ \hline
$ 2^{44}$= 17592186044416&  3.88364187E-0009 &  9.05905329E-0007 &  4.28702840E-0003 \\ \hline
\end{tabular}
\end{center}

It's clear that the error term ratio stays fairly constant, so that the order of magnitude is correct, although the numerical constants in the error formula need considerable improvement!
\subsubsection*{Acknowledgment}
I thank Joseph C. V\'arilly for comments on an earlier version. 
Support from the Vicerrector\'{\i}a de Investigaci\'on of the 
University of Costa Rica is acknowledged.


\end{document}